\documentclass[12pt, equation, leqno]{article}
\usepackage{amssymb}
\usepackage{amsmath}
\usepackage{theorem}
\evensidemargin\oddsidemargin
\makeatletter
\@addtoreset{equation}{section}

\makeatother
\newtheorem{theo}{Theorem}

\newtheorem{lem}{Lemma}
\newtheorem{cor}{Corollary}
{\theorembodyfont{\rmfamily}  }
\newcommand{\PP}{\Bbb{P}}

\newcommand{\E}{\Bbb{E}}

\newcommand{\be}{\begin{equation}}
\newcommand{\ee}{\end{equation}}
\newcommand{\nin}{\noindent}
\begin{document}

\begin{center}
{\bf CHARACTERIZATION OF LIL BEHAVIOR IN BANACH SPACE}

\vskip 0.3cm

\renewcommand{\thefootnote}{\fnsymbol{footnote}}
UWE EINMAHL$^{a,}$\footnote[1]{Research supported  by an FWO Grant} and DELI LI$^{b,}$\footnote[2]{Research supported by a grant from 
the Natural Sciences and Engineering Research Council of Canada. }

\vskip 0.3cm

$^a$ {\it Departement  Wiskunde, Vrije Universiteit Brussel,} \\
{\it Pleinlaan 2, B-1050 Brussel, Belgium;} \\
$^b$ {\it Department of Mathematical Sciences, Lakehead University,}\\
{\it Thunder Bay, Ontario, Canada P7B 5E1}
\end{center}

\begin{abstract}
\nin
In a recent paper by the authors a general result characterizing
two-sided LIL behavior for real valued random variables has been
established. In this paper, we look at the corresponding problem in the Banach space setting. We show that there are analogous results in this more general setting. In particularly, 
we provide a necessary and sufficient condition for 
LIL behavior with respect to  sequences of the form
$\sqrt{nh(n)}$, where $h$ is from a suitable subclass of the
positive, nondecreasing slowly varying functions. To prove these results we have to use a different method.  
One of our main tools is an improved Fuk-Nagaev type inequality in Banach space which should be of independent interest.
\end{abstract}

\bigskip\noindent {\it Short title:} LIL behavior in Banach Spaces



\newpage

\section{Introduction}

Let $(B, \| \cdot \| )$ be a real separable Banach space
with topological dual $B^{*}$. Let $ \{X,~X_{n};~n \geq 1 \}$ 
be a sequence of independent and identically distributed (i.i.d.)
B-valued random variables. As usual, let   
$S_{n} = \sum_{i=1}^{n} X_{i},~ n \geq 1$ and
set $Lt=\log(t\vee e), LLt = L(Lt), t \ge 0.$

One of the classical results of probability is the Hartman-Wintner LIL and the definitive version of this result in Banach space has been proven by Ledoux and Talagrand (1988). \vspace{.2cm}

\nin{\bf Theorem A} {\it A random variable $X: \Omega \to B$ satisfies the bounded LIL, that is
\be \limsup_{n \to \infty} \|S_n\|/\sqrt{nLLn} < \infty \mbox{ a.s.}\label{LT0}\ee
if and only if the following three conditions are fulfilled:
\be \E\|X\|^2/LL\|X\| < \infty, \;\;\E X =0, \label{LT1}\ee
\be \E f^2(X) < \infty, f \in B^*, \label{LT2}\ee
\be \{S_n/\sqrt{nLLn}\} \mbox{ is bounded in probability}.\label{LT3}\ee}
Furthermore it is known that if one assumes instead of (\ref{LT3}),
\be S_n/\sqrt{nLLn} \stackrel{\PP}{\to} 0 \label{LT4}\ee
one has
\be \limsup_{n \to \infty} \|S_n\|/\sqrt{2nLLn} = \sigma \mbox{ a.s.} \label{LT5}\ee
where $\sigma^2 = \sup_{f \in B_1^*} \E f^2(X)$ and $B_1^*$ is the unit ball of $B^*.$ It is easy to see that $\sigma^2$ is finite under assumption (\ref{LT2}).\\
If $B$ is a type 2 space then (\ref{LT1}) implies (\ref{LT4}) and the bounded LIL holds if and only if conditions (\ref{LT1}) and (\ref{LT2}) are satisfied. Moreover, in this case we also know  the exact value of the limsup in (\ref{LT0}). \\
Recall that we call a Banach space type 2 space  if we have for any sequence $\{Y_n\}$ of independent mean zero random variables with $\E\|Y_n\|^2 < \infty, n \ge 1:$
$$\E \|\sum_{i=1}^n Y_i\|^2 \le C \sum_{i=1}^n \E \|Y_i\|^2, n \ge 2.$$
where $C > 0$ is a constant. It is well known that finite-dimensional spaces and Hilbert spaces are type 2 spaces. \vspace{.3cm}

Finding the precise value of $\limsup_{n \to \infty}\|S_n\|/\sqrt{2nLLn}$ in general seems to  be a difficult problem (see, for instance, Problem 5 on page 457 of Ledoux and Talagrand (1991)). If one imposes the stronger assumption $\E \|X\|^2 < \infty$ and $\E X=0$ instead of (\ref{LT1})  and (\ref{LT2}), de Acosta, Kuelbs and Ledoux (1983) proved that with probability one,
\be \sigma \vee \beta_0 \le \limsup_{n \to \infty}\|S_n\|/\sqrt{2nLLn} \le \sigma + \beta_0, \label{alt} \ee
where $\beta_0 = \limsup_{n \to \infty} \E\|S_n\|/\sqrt{2nLLn}.$ Moreover, they showed that the lower bound $\sigma \vee \beta_0$ is sharp for random variables in $c_0.$ It is still open whether this is the case in other Banach spaces as well.
If $S_n/\sqrt{nLLn} \stackrel{\PP}{\to} 0$, one has $\beta_0 = 0$ and one can re-obtain result (\ref{LT5}) if $\E \|X\|^2 < \infty$. Also note that in all other cases one misses the ``true'' value of the $\limsup$ at most by  a factor 2. So if $\E \|X\|^2 < \infty,$ we have  a fairly complete picture and it is natural to ask whether it is possible to establish (\ref{alt})  under   conditions (\ref{LT1}) and (\ref{LT2}). This has  been shown by de Acosta, Kuelbs and Ledoux (1983) for certain Banach spaces which satisfy a so-called upper Gaussian comparison principle, but the question of whether this is the case for general Banach spaces seems to be still open. As a by-product of our present work  we will  be able to answer this  in the affirmative.\\
 
There are also extensions of the Hartman-Wintner LIL to real-valued random variables with possibly infinite variance. Feller (1968) obtained an LIL for certain variables in the domain of attraction to the normal distribution and this was further generalized by Klass (1976, 1977).  Kuelbs (1985) and Einmahl (1993) found versions of these results in the Banach space setting. In a recent paper Einmahl and Li (2005) looked at the the following problem for real-valued random variables:\vspace{.2cm}

\nin \emph {Given a sequence, $a_n=\sqrt{nh(n)},$ where $h$ is a slowly varying non-decreasing function, when does one have with probability one,}
$0 < \limsup_{n \to \infty}|S_n|/a_n < \infty ? $ \vspace{.2cm}

Somewhat unexpectedly it turned out that the classical Hartman-Wintner LIL could be generalized to a ``law of the very slowly varying function''. It is the main purpose of the present paper to investigate whether there are also such results in the Banach space setting. In the process we will derive a very general result on almost sure convergence (see Theorem 5, Sect. 3) which specialized to the classical normalizing sequence $\sqrt{2nLLn}$ also gives result (\ref{alt}) under the weakest possible conditions.
\section{Statement of main results}
Let ${\cal H}$ be the set of all continuous, non-decreasing 
functions $h: [0,\infty) \to (0,\infty)$, which are slowly
varying at infinity. To simplify notation we set  $\Psi(x)=  \sqrt{xh(x)}, x \ge 0$ for $h \in \mathcal{H}$ and let $a_n = \Psi(n), n \ge 1.$\\
Given a random variable $X: \Omega \to B$ we consider an infinite-dimensional truncated second moment function $H: [0,\infty) \to [0,\infty)$ defined by 
$$ H(t) :=\sup_{f \in B_1^*} \E f^2(X)I\{\|X\| \le t\}, t \ge 0.$$
The first theorem gives a characterization for having  $\limsup_{n \to \infty}\|S_n\|/a_n < \infty$ a.s. where $a_n$ is a normalizing sequence of the above form.\vskip 0.3cm

\begin{theo}
Let $ X $ be a B-valued random variable.
Then we have
\begin{equation}\label{2.1}
\limsup_{n \rightarrow \infty}
\frac{\|S_{n}\|}{a_{n}} < \infty ~~\mbox{a.s.}
\end{equation}
if and only if 
\begin{equation}\label{2.2}
\E X = 0, ~~~\E \Psi^{-1}(\|X\|) < \infty, 
\end{equation}
\begin{equation}\label{2.3}
\mbox{the sequence} ~~
\{S_{n}/a_{n};~n \geq 1 \}
~\mbox{is bounded in probability,}
\end{equation}
and there exists $ c \in [0, \infty)$ such that  
\begin{equation}\label{2.4}
\sum_{n=1}^{\infty} \frac{1}{n} \exp \left \{ 
- \frac{c^{2} h(n)}{2H(a_{n})} \right \} < \infty.
\end{equation}
\end{theo}

\nin By strengthening  condition (\ref{2.3}) we can find the exact limsup value in (\ref{2.1}).
\begin{theo}
Assume (\ref{2.2}) holds and (\ref{2.3}) is strengthened to
\begin{equation}\label{2.5}
S_{n}/a_{n} \stackrel{\PP}{\to} 0,
\end{equation}
then
\begin{equation}\label{2.6}
\limsup_{n \rightarrow \infty}
\frac{\|S_{n}\|}{a_{n}} = C_{0} ~~\mbox{a.s.},
\end{equation}
where 
\begin{equation}\label{2.7}
C_{0} = \inf \left \{c \geq 0:~\sum_{n=1}^{\infty} \frac{1}{n} \exp \left \{
- \frac{c^{2} h(n)}{2H(a_{n})} \right \} < \infty \right \}.
\end{equation}
\end{theo}  
As in the classical case (when considering the sequence $a_n = \sqrt{2nLLn}$ ) one can show that in type 2 spaces (\ref{2.2}) implies (\ref{2.5}) so that in this case 
(\ref{2.1}) holds if and only if conditions (\ref{2.2}) and (\ref{2.4}) are satisfied. Moreover, the value of the $\limsup$ in (\ref{2.1}) is then always equal to $C_0$. \\

In general, it can be difficult to determine this parameter. For this reason 
we now look at normalizing sequences $a_n=\sqrt{nh(n)}$ for functions $h$ from  certain subclasses of $\mathcal{H}.$ Given $0 \le q < 1$, 
let  ${\cal H}_q \subset  {\cal H}$ the class
which contains all functions $h \in \mathcal{H}$ satisfying the condition,
\[
\lim_{t \to \infty} \frac{h(tf_{\tau}(t))}{h(t)} =1,
~~0 < \tau < 1 - q,
\] 
where $f_{\tau} (t) = \exp ((Lt)^{\tau}), ~0 \le \tau\le 1$. Finally let ${\cal H}_1 = {\cal H}$.
Clearly ${\cal H}_{q_{1}} \subset {\cal H}_{q_{2}}$ whenever
$0 \leq q_{1} < q_{2} \leq 1$. We call the functions in the smallest subclass ${\cal H}_0$ ``very slowly varying''.  
From the following theorem it follows that under assumption (\ref{2.5}) we  have $C_0 \le \lambda$ for {\it any} $h \in \mathcal{H}$ where  $\lambda$ is a parameter which can be easily determined via the  $H$-function. If we have $h \in \mathcal{H}_q$, then it also follows that $C_0 \ge (1-q)^{1/2}\lambda.$ Thus, if $h \in \mathcal{H}_0,$ we have $C_0=\lambda$ and this way we can extend the classical LIL to a ``law of the very slowly varying function''.  Possible choices for very slowly varying functions are for instance  $(LLt)^p,\, p \ge 1$ and $(Lt)^r,\, r >0.$
\begin{theo}
Let $X$ be a B-valued random variable.
Suppose now that $h \in {\cal H}_q$ where
$0 \le q \le 1.$ Assume (\ref{2.2}) and (\ref{2.5}) hold.
Then 
\begin{equation}\label{2.8}
(1 - q)^{1/2} \lambda \leq 
\limsup_{n \rightarrow \infty}
\frac{\|S_{n}\|}{a_{n}} \leq \lambda ~~\mbox{a.s.,}
\end{equation}
where  
\begin{equation}\label{2.9}
\lambda^{2} = \limsup_{x \rightarrow \infty}
\frac{2\Psi^{-1}(x LL x)}{x^{2} LLx}H(x). 
\end{equation}
\end{theo}
Note the $\limsup$ in condition ({\ref{2.9}). If this $\limsup$ is actually a limit, then it easily follows from Theorem 2 that  $\limsup_{n\to \infty}\|S_n\|/a_n = \lambda$ a.s. for any function $h \in \mathcal{H}$. This condition, however, is  not necessary.   It is a special feature of the function class $\mathcal{H}_0$ that under condition (\ref{2.5}) the  $\limsup$ in (\ref{2.9}) being equal to $\lambda^2$ is  necessary and also sufficient in combination with (\ref{2.2})  for having $\limsup_{n\to \infty}\|S_n\|/a_n = \lambda$ a.s .  Moreover, we have for $h \in \mathcal{H}_q$ and $0 \le q <1$ that  $\limsup_{n\to \infty}\|S_n\|/a_n < \infty$ a.s. if and only if $\lambda < \infty$ and conditions (\ref{2.2}) and (\ref{2.3}) hold.\\ 

 Theorem 3 gives us  analogous corollaries as in the real-valued case. We state two of these. The formulation of the other ones, for instance,  a law of the logarithm (see, Corollary 2, Einmahl and Li (2005)), should be then obvious.
\begin{cor}
Let $X$ be a B-valued random variable.
Let $p \geq 1$. Then we have
\begin{equation}\label{2.10}
\limsup_{n \rightarrow \infty}
\frac{\|S_{n}\|}{\sqrt{2n (LLn)^{p}}}
< \infty ~~\mbox{a.s.}
\end{equation}
if and only if 
\begin{equation}\label{2.11}
\E X = 0, ~~~\E \|X\|^{2}/(LL\|X\|)^{p} < \infty,
\end{equation}
\begin{equation}\label{2.12}
\lambda^{2} = \limsup_{x \rightarrow \infty}
(LLx)^{1-p} H(x) < \infty,
\end{equation}
and
\begin{equation}\label{2.13}
\mbox{the sequence} ~~
\{S_{n}/\sqrt{2n (LLn)^{p}};~n \geq 1 \}
~\mbox{is bounded in probability.}
\end{equation}
Furthermore,
\begin{equation}\label{2.14}
\limsup_{n \rightarrow \infty}
\frac{\|S_{n}\|}{\sqrt{2n (LLn)^{p}}}
= \lambda~~\mbox{a.s.}
\end{equation}
whenever  condition (\ref{2.14}) is strengthened to
\begin{equation}\label{2.15}
S_{n}/\sqrt{2n (LLn)^{p}} \stackrel{\PP}{\to} 0.
\end{equation}
\end{cor} 
If $p=1$ we re-obtain Theorem A, but the above corollary actually shows that we have for any $p \ge 1$ an LIL. 
If $\lambda = 0$ in Theorem 3, we  obtain the following useful stability result. 
\begin{cor}
Assume that $X: \Omega \to B$ is a random variable satisfying
\be \E X = 0, ~~~\E \Psi^{-1}(\|X\|) < \infty, \label{s1}\ee
\be \lim_{x \rightarrow \infty}
\frac{\Psi^{-1}(x LL x)}{x^{2} LLx} H(x) = 0, \label{s2}\ee
\be S_n/a_n \stackrel{\PP}{\to} 0 \label{s3}\ee
then we have
\begin{equation}\label{s4}
\lim_{n \rightarrow \infty}
\frac{S_{n}}{a_{n}} = 0 ~~\mbox{a.s.}
\end{equation}
Conversely, if $q < 1$ then (\ref{s4}) implies (\ref{s1}) - (\ref{s3}).
\end{cor}
The remaining part of the paper is organized as follows. In Sect. 3 we state and prove an infinite-dimensional version of the Fuk-Nagaev inequality improving an earlier version of this inequality given as Theorem 5 in Einmahl (1993). Using a recent result of Klein and Rio (2005) who obtained in some sense an optimal version of the classical Bernstein inequality in infinite-dimensional spaces,  we can replace the constant 144 in the exponential term of the earlier version by  $2+\delta$ for any $\delta >0$. Employing this improved version of the Fuk-Nagaev inequality one can give much more direct proofs for LIL results than in Einmahl (1993). Especially it is no longer necessary to use randomization arguments and Sudakov minoration  for obtaining the precise value of $\limsup_{n \to \infty}\|S_n\|/a_n.$ Readers who are mainly interested in inequalites can read this part independently of the other parts of the present paper. In Sect. 4 we then use the improved Fuk-Nagaev inequality to establish the upper bound part of a general result on  almost sure convergence for normalized sums $S_n/c_n$ where $\{c_n; n \ge 1\}$ is a sufficiently regular normalizing sequence. This includes all sequences $a_n =\sqrt{nh(n)}$, where $h \in \mathcal{H}.$  For proving the lower bound part we first use  an extension of a method employed in the proof of Theorem 2, Einmahl (1993) to get a first lower bound (see Section 4.2). In the classical case $c_n=\sqrt{2nLLn}$ this bound would be equal to $\sigma$. Our method is fairly elementary and one only needs classical results such as a non-uniform bound on the convergence speed for the CLT on the real line. In Sect. 4.3 we obtain a second lower bound which, in the classical case, matches $\beta_0$ defined in (\ref{alt}). Here we use a modification of an argument based on Fatou's lemma which is due to de Acosta, Kuelbs and Ledoux (1983).  In Sect. 5 we finally infer the results stated in Sect. 2 from our general almost sure convergence result (Theorem 5).

\section{A Fuk-Nagaev type inequality}
As mentioned in Sect. 2 we use an infinite-dimensional version of the Bernstein inequality which essentially goes back to Talagrand (1994). This inequality turned out to be extremely useful in many applications, but there was a shortcoming that there were no explicit numerical constants.  Ledoux (1996) found a different and very elegant method for proving such inequalities  which is based on a log-Sobolev type argument in combination  with a tensorization of the entropy. He was also able to provide concrete numerical constants for these inequalites. His method was subsequently refined by Massart (2000) and Rio (2002)  among other authors.  Bousquet (2002)  obtained optimal constants in the iid case.  Finally, Klein and Rio (2005) generalized this result to independent, not necessarily identically distributed random variables. Their results are formulated for empirical processes, but using a standard argument one can easily obtain  inequalities for sums of independent $B$-valued variables from the ones for empirical processes. \vspace{.2cm}

\nin We need the following fact which follows  from  Lemma 3.4 of Klein and Rio (2005).  \\
{\bf Fact A } {\it Let $Y_1,\ldots, Y_n$ be independent  $B$-valued random variables with mean zero such that 
$$\|Y_i\| \le M \mbox{ a.s.},1 \le i \le n. $$
 Then we have for $0 < s < 2/(3M)$:
 \be
 \E \exp(s\|\sum_{i=1}^n Y_i\|) \le \exp\left(s\E \|\sum_{i=1}^n Y_i\| + \beta_n s^2/(2 - 3Ms)\right) \label{KR}
\ee
where $\beta_n = 2M\E\|\sum_{i=1}^n Y_i\| + \Lambda_n$ with $\Lambda_n = \sup\{\sum_{j=1}^n \E f^2(Y_j) : f \in B_1^*\}$ and $B_1^*$ is equal to the unit ball of} $B^*$.\vspace{.3cm}

To prove this inequality we set $Z_i =Y_i/M, 1\le i \le n$.
Recall that $B$ is separable so that we have for any $z \in B$,  $\|z\| = \sup_{f \in D} f(z),$  where $D$ is a countable subset of $B_1^*.$ Set in Theorem 1.1 of Klein and Rio (2005) $\mathcal{X}=B$ and consider the following countable class of functions from $\mathcal{X}$ to $[-1,1]^n$:
$\mathcal{S} = \{(-1\vee(f \wedge 1),\ldots, -1\vee(f \wedge 1)): f \in D\}.$
Then we readily obtain that $\sup_{s \in \mathcal{S}} \{s^1(Z_1) + \ldots +s^n(Z_n)\} =\|Z_1 + \ldots +Z_n\|$ a.s.  and we can infer from the afore-mentioned lemma that for $0 < t <2/3,$
$\E \exp(t\|\sum_{i=1}^n Z_i\|) \le \exp\left(t\E \|\sum_{i=1}^n Z_i\| + \gamma_n t^2/(2 - 3t)\right),$
where $\gamma_n = 2\E \|\sum_{i=1}Z_i\| + V_n$ and $V_n = \sup\{\sum_{j=1}^n \E f^2(Z_j) : f \in B_1^*\}.$ Replacing $Z_i$ by $Y_i/M$ and setting $s=t/M$ we obtain (\ref{KR}).\vspace{.2cm}

\nin Using the well known fact that $\exp(s\|\sum_{i=1}^k Y_i\|), 1 \le k \le n$ is a submartingale if $s >0$ (recall that we assume $\E Y_k = 0, 1 \le k \le n$), we can infer from Doob's maximal inequality for submartingales that for any $x >0,$
$$\PP\left\{\max_{1 \le k \le n}\|\sum_{i=1}^k Y_i\| \ge \E \|\sum_{i=1}^n Y_i\| + x\right\} \le \exp( \beta_n s^2/(2 - 3Ms) - sx), 0 < s < 2/(3M).$$
Choosing $s = 2x/(2\beta_n + 3Mx)$ we finally obtain that
\be \PP\left\{\max_{1 \le k \le n}\|\sum_{i=1}^k Y_i\| \ge \E \|\sum_{i=1}^n Y_i\| + x\right\} \le
\exp\left(-\frac{x^2}{2\Lambda_n + (4\E\|\sum_{i=1}^n Y_i\| + 3x)M}\right).\label{KR1}\ee
Next note that we trivially have for any $\epsilon >0,$
\begin{eqnarray} \label{triv}
&&\exp\left(-\frac{x^2}{2\Lambda_n + (4\E\|\sum_{i=1}^n Y_i \| +3x)M}\right)\\
&\le& \exp\left(-\frac{x^2}{(2+\epsilon)\Lambda_n} \right) +
\exp\left(-\frac{x^2}{(1+2/\epsilon) (4\E\|\sum_{i=1}^n Y_i \| +3x)M}\right).\nonumber
\end{eqnarray}
Combining (\ref{KR1}) and (\ref{triv}) and setting $x= \eta \E\|\sum_{i=1}^n Y_i \| +  y,$
where $0 <\eta \le 1$ and $y>0,$ we can conclude that for any $y >0,$
\be
\PP\left\{\max_{1 \le k \le n}\|\sum_{i=1}^k Y_i \| \ge (1 + \eta)\E\|\sum_{i=1}^n Y_i \| + y\right\}
\le \exp\left(-\frac{y^2}{(2+\epsilon)\Lambda_n} \right) +
\exp\left(-\frac{y}{D_{\epsilon,\eta}\,M}\right), \label{KRmod}
\ee
where $D_{\epsilon,\eta} = (1+2/\epsilon)(3 + 4/\eta).$
We are now ready to prove 
\begin{theo} \label{FNa}
Let $Z_1,\ldots, Z_n$ be independent B-valued random variables with mean zero such that for some $s >2,$ $\E\|Z_i\|^s < \infty, 1 \le i \le n.$ Then we have for $0 < \eta \le 1, \delta >0$ and any $t >0,$
\be
\PP\left\{\max_{1 \le k \le n}\|\sum_{i=1}^k Z_i \| \ge (1+\eta)\E\|\sum_{i=1}^n Z_i \| + t\right\}
\le \exp\left(-\frac{t^2}{(2+\delta)\Lambda_n} \right) + C\sum_{i=1}^n \E\|Z_i\|^s/t^s,
\ee
where $\Lambda_n = \sup\{\sum_{j=1}^n \E f^2(Z_j) : f \in B_1^*\}$ and $C$ is a positive constant depending on $\eta,\delta$ and $s$.
\end{theo}
{\bf Proof.} To simplify notation we set for $y >0$
$$\beta(y)=\beta_s (y)  = \sum_{i=1}^n \E \|Z_i\|^s/y^s.$$
Assume that $\beta(y) < 1.$ For $\epsilon > 0$  fixed we consider the following truncated variables 
$$Y_i := Z_i I\{ \|Z_i\| \le \rho\epsilon y\}, \,Y'_i = Y_i - \E Y_i,  1 \le i \le n,$$
 where 
$$\rho= \rho(\epsilon,\eta, y)= 1 \wedge \frac{1}{2\epsilon D_{\epsilon,\eta}\,\log(1/\beta(y))} .$$
Applying inequality (\ref{KRmod}) with $M=2\rho \epsilon y$ we find that
\be
\PP\left\{\max_{1 \le k \le n}\|\sum_{i=1}^k Y'_i \| \ge (1+\eta)\E\|\sum_{i=1}^n Y'_i \| + y\right\}
\le \exp\left(-\frac{y^2}{(2+\epsilon)\Lambda_n} \right) + \beta(y). \label{in1}
\ee
Next consider the variables
$$\Delta_i := Z_i I\{\rho\epsilon y < \|Z_i\| \le \epsilon y\}, 1 \le i \le n.$$
Employing the Hoffmann-J\o rgensen inequality (see, for instance, inequality (6.6) in Ledoux and Talagrand (1991)), we can conclude that
\be
\PP\left\{ \max_{1 \le k \le n}\|\sum_{i=1}^k \Delta_i\| \ge 4\epsilon y \right\}
\le \left(\PP\left\{\max_{1 \le k \le n} \|\sum_{i=1}^k \Delta_i\| \ge \epsilon y \right\}\right)^2 \label{in2}
\ee
which in turn is 
$$\le \left(\sum_{i=1}^n \PP\{\Delta_i \ne 0\} \right)^2
\le \left(\sum_{i=1}^n \PP\{\|Z_i\| \ge  \rho \epsilon y\} \right)^2.$$
Using Markov's inequality and recalling the definition of $\rho$ we see that this last term is bounded above by
$$
(2D_{\epsilon,\eta})^{2s}\beta^2(y)(\log(1/\beta(y)))^{2s}
\le K_s (2D_{\epsilon,\eta})^{2s}\beta(y),
$$
where $K_s >0$ is a constant so that $(\log a)^{2s}\le K_s a , a \ge 1.$
We can conclude that
\be
\PP\left\{\max_{1 \le k \le n} \|\sum_{i=1}^k \Delta_i\| \ge 4\epsilon y \right\} \le C' \beta(y),
\label{in3}
\ee
where $C'=K_s(2D_{\epsilon,\eta})^{2s}.$ \\
Next set $\Delta'_i := Z_i I\{\|Z_i\| > \epsilon y\}, 1 \le i \le n.$ Then we have once more by Markov's inequality
\be
\PP\left\{\max_{1 \le k \le n}\|\sum_{i=1}^k \Delta'_i\| \ne 0\right\}\le \epsilon^{-s}\beta(y).
\label{in4}
\ee
Combining inequalities (\ref{in1}), (\ref{in3}) and (\ref{in4}), we see that if $\beta(y) < 1$ we have
$$
\PP\left\{\max_{1 \le k \le n} \|\sum_{i=1}^k (Z_i - \E Y_i) \| \ge (1+\eta)\E\|\sum_{i=1}^n Y'_i \| + (1+4\epsilon)y\right\}
\le \exp\left(-\frac{y^2}{(2+\epsilon)\Lambda_n} \right) + C''\beta(y), $$
where $C''=1+ C' + \epsilon^{-s}$. 
A simple application of the triangular inequality gives 
\be
\PP\left\{\max_{1 \le k \le n} \|\sum_{i=1}^k Z_i  \| \ge b'_n+ (1+4\epsilon)y\right\}
\le \exp\left(-\frac{y^2}{(2+\epsilon)\Lambda_n} \right) + C''\beta(y), \label{in5}
\ee
where \vspace{-1cm}

\begin{eqnarray*}
b'_n &=& (1+\eta)\E\|\sum_{i=1}^n Y'_i\| + \max_{1 \le k \le n} \|\sum_{i=1}^k \E Y_i\|\\
& \le&
(1+\eta)\E\|\sum_{i=1}^n Y_i\| + 3\max_{1 \le k \le n} \|\sum_{i=1}^k \E Y_i\|.
\end{eqnarray*}
Further note that \vspace{-.5cm}

\begin{eqnarray*}
\E\|\sum_{i=1}^n Y_i \| &\le & \E \|\sum_{i=1}^n Z_i \| + \E\|\sum_{i=1}^n Z_i I\{\|Z_i\| \ge  \rho\epsilon y\} \|\\
&\le&\E \|\sum_{i=1}^n Z_i \| + \sum_{i=1}^n \E\|Z_i\|I\{\|Z_i\| \ge \rho \epsilon y\}
=: \E \|\sum_{i=1}^n Z_i \| + \delta_n.
\end{eqnarray*}
As we have $\E Z_i =0, 1 \le i \le n$ it also follows that $\max_{1 \le k \le n} \|\sum_{i=1}^k \E Y_i\| \le \delta_n$ and consequently,
\be b'_n \le (1+\eta)\E\|\sum_{i=1}^n Z_i\| + 5\delta_n. \label{cent}\ee
Furthermore, we have,
$$ \delta_n \le y\beta(y)/\{\rho \epsilon\}^{s-1} \le \epsilon y$$
provided that $\beta(y) \le \epsilon^s \rho^{s-1}.$\\ It is easily checked that if $\rho <1,$ we have  $ \beta(y)/(\epsilon^s \rho^{s-1}) \le \beta(y)/(\epsilon^s \rho^{s}) \le (C''\beta(y))^{1/2}.$ (We are assuming that $\beta(y) \le 1.$) Consequently,    $\delta_n \le \epsilon y$ whenever $C''\beta(y) \le 1$ and $\rho <1$. This is also true if $\rho =1$ as we have $C'' \ge \epsilon^{-s}.$ We thus can conclude if $\beta(y) \le 1/C'' < 1:$ 
\be
\PP\left\{\max_{1 \le k \le n}\|\sum_{i=1}^k Z_i \| \ge (1+\eta)\E\|\sum_{i=1}^n Z_i \| + (1+9\epsilon)y\right\}
\le \exp\left(-\frac{y^2}{(2+\epsilon)\Lambda_n} \right) + C''\beta(y). 
\label{in6}
\ee
The above inequality is of course trivial if $\beta(y) > 1/C''$ and  consequently (\ref{in6}) holds for all $y >0.$ 
Setting $y= t/(1+9\epsilon)$ and choosing $\epsilon$ in (\ref{in6}) so small that $(2+\epsilon)(1+9\epsilon)^2 \le 2+\delta$, we obtain the assertion. $\Box$
\section{A general  result on almost sure convergence}
 Let  $c_n$ be a sequence of real numbers satisfying the following two conditions,
\be
c_n/\sqrt{n} \nearrow \infty \label{RE}
\ee
and
\be
 \forall\; \epsilon >0\,\exists\, m_\epsilon \ge 1: c_n/c_m \le (1+\epsilon)(n/m), m_\epsilon \le m < n. \label{REG}
\ee  
Note that condition (\ref{REG}) is satisfied for any sequence $a_n$ considered in Section 2.  \\
Let $H$ be defined as in  Section 1, that is 
$$H(t) = \sup_{f \in B_1^*}\E f^2(X) I\{\|X\| \le t\}, t \ge 0.$$
Set 
$$\alpha_0 = \sup\left\{\alpha \ge 0: \sum_{n=1}^{\infty} n^{-1}\exp\left(-\frac{\alpha^2 c^2_n}{2nH(c_n)}\right) = \infty\right\}
.$$
In general $\alpha_0$ can be any number in $[0,\infty].$ If we are assuming that $\E f^2(X) < \infty, f \in B^*$ and we choose $c_n = \sqrt{2nLLn}$, it follows that $\alpha_0^2 = \sigma^2 = \sup_{f \in B^*_1}\E f^2(X).$\\

Our main result in this section is the following generalization of (\ref{alt}),
\begin{theo}\label{th1} Let $X, X_1, X_2, \ldots$ be i.i.d. mean zero random variables taking values in a separable Banach space $B$. 
Assume that 
\be
\sum_{n=1}^{\infty} \PP\{\|X\| \ge c_n\} < \infty, \label{MOM}
\ee
where $c_n$ is a  sequence of positive real numbers satisfying conditions (\ref{RE}) and (\ref{REG}).\\
Then we have with probability one,
\be
\alpha_0 \vee \beta_0  \le \limsup_{n \to \infty}\|S_n\|/c_n \le \alpha_0 + \beta_0,\label{bounds}
\ee
where $\beta_0 = \limsup_{n \to \infty} \E \|S_n\|/c_n.$ 
\end{theo}
The following lemma which is  more or less known shows that $\beta_0$ is finite whenever $\{S_n/c_n; n \ge 1\}$ is bounded in probability and that $\beta_0 =0$ if $S_n/c_n \stackrel{\PP}{\to} 0.$ So in the latter case we  see that the $\limsup$ in (\ref{bounds}) is equal to $\alpha_0.$ 
\begin{lem} Let $X, X_1, X_2, \ldots$ be iid B-valued random variables with mean zero and let $S_n = \sum_{i=1}^n X_i, n \ge 1$.  Let $\{c_n\}$ be a sequence of positive  real numbers satisfying conditions (\ref{RE}) and (\ref{REG}). Under assumption (\ref{MOM}) we have the following equivalences:
\begin{itemize}
\item[(a)] $\{S_n/c_n; n \ge 1\}$ is bounded in probability $\Longleftrightarrow \limsup_{n \to \infty} \E\|S_n\|/c_n < \infty.$ 
\item[(b)] $S_n/c_n \stackrel{\PP}{\to} 0 \Longleftrightarrow \E\|S_n\|/c_n \to 0.$
\end{itemize}
\end{lem}
{\bf Proof} We only need to prove the implications ``$\Rightarrow$'' and by a standard symmetrization argument it is enough to do that for symmetric random variables.
We have for any $\epsilon >0,$
\be
\E \|S_n\| \le \E\left\|\sum_{i=1}^n X_i I\{\|X_i\| \le \epsilon c_n\}\right\| +
n \E \|X\|I\{\|X\| > \epsilon c_n\}. \label{last}
\ee
The last term is of order $o(c_n)$ under assumption (\ref{MOM}) (see Lemma 1, Einmahl and Li (2005)). 
Using the trivial inequality
$$\PP\left \{\|\sum_{i=1}^n X_i I\{\|X_i\| \le \epsilon c_n\}\| \ge x\right\} 
\le \PP\{\|S_n\| \ge x\} + n \PP\{\|X\| \ge \epsilon c_n\},$$
in conjunction with Proposition 6.8 in Ledoux and Talagrand (1991), we find that if $\{S_n/c_n\}$ is bounded in probability, the first term in (\ref{last}) is $\le C(\epsilon) < \infty$ . Consequently, we have in this case, $\E\|S_n\|/c_n < \infty.$  Assuming $S_n/c_n \stackrel{\PP}{\to} 0,$  one can choose $C(\epsilon)$ so that  $C(\epsilon) \to 0$ as $\epsilon \to 0$. Since we can make $\epsilon$ arbitrarily small, it follows that $\E\|S_n\|/c_n \to 0$ if $S_n/c_n \stackrel{\PP}{\to} 0. \;\Box$\\

 If $B$ is a type 2 Banach space, assumption (\ref{MOM}) implies that $\E \|S_n\| = o(c_n).$ (See Lemma 6, Einmahl (1993). The proof given there works also under the present conditions on $\{c_n\}.$) Therefore we have in any type 2 space, $\beta_0=0$ and  the $\limsup$ in (\ref{bounds}) is equal to $\alpha_0$. Recalling that finite dimensional spaces are type 2 spaces, we see that this result extends Theorem 3 of Einmahl and Li (2005). Also note that the conditions on $\{c_n\}$ are general enough so that one  can  infer Theorem 3, Einmahl (1993) from the present Theorem 6 as well (without using randomization and Sudakov minoration). \vspace{.3cm}
 
We now turn to the proof of Theorem \ref{th1}. We assume throughout that condition (\ref{MOM}) is satisfied.
Using essentially the same argument as in Lemma 3 of Einmahl (2007) one can infer from the definition of $\alpha_0$ that whenever
$n_j \nearrow \infty$ is a subsequence satisfying for large enough $j,$
\be 1 < a_1 < n_{j+1}/n_j \le a_2 < \infty,\label{sub}\ee
we have:
\be
\sum_{j=1}^{\infty} \exp\left(-\frac{\alpha^2 c_{n_j}^2}{2n_j H(c_{n_j})}\right)  \begin{cases} =\infty & \mathrm{if}\;\alpha < \alpha_0,\\ <\infty & \mathrm{if}\;\alpha > \alpha_0.\end{cases} \label{geom}
\ee

\subsection{The upper bound part}
W.l.o.g. we can and do assume in this part that $\alpha_0 + \beta_0 < \infty.$\\ 
We first note that on account of (\ref{RE})   and assumption (\ref{MOM}) we have   for any subsequence $\{n_j\}$ satisfying  (\ref{sub}),
\be
\sum_{j=1}^{\infty} n_j \E \|X\|^3I\{\|X\| \le c_{n_j}\}/c_{n_j}^3 < \infty.
\label{fact1}
\ee 
(See, for instance, Lemma 7.1, Pruitt (1981).)\\
Moreover, we have as $n \to \infty,$
\be
n\E \|X\|I\{\|X\| \ge c_n\} \le \sum_{i=1}^n \E\|X\|I\{\|X\| \ge c_i\} = o(c_n)
\label{fact2}
\ee
 This last fact follows as in the proof of Lemma 10, Einmahl (1993) replacing $\gamma_n$ by $c_n$. \\
Set $X'_n = X_n I\{\|X_n\| \le c_n\}, n \ge 1$ and denote the sum of the first $n$ of these variables by $S'_n, n \ge 1.$
We obviously have 
$$\sum_{n=1}^{\infty}\PP\{X_n \ne X'_n\} < \infty$$
so that with probability one, $X_n = X'_n$ eventually. Due to relation (\ref{fact2}) we have $\E S'_n = o(c_n)$ and consequently it is enough to show,
\be
\limsup_{n \to \infty}\|S'_n - \E S'_n\|/c_n \le \alpha_0 + \beta_0 \mbox{ a.s.}\label{upper}
\ee
This follows via Borel-Cantelli once we have proven for any $0 <\delta <1$
\be
\sum_{j=1}^{\infty}\PP\left\{\max_{1 \le n \le n_{j+1}}\|S'_n -\E S'_n\| \ge \{\alpha_0  +\delta +\beta_0(1 +\delta)\}(1+ 2\delta)c_{n_j}\right\} < \infty, \label{series}
\ee
where $n_j \sim \rho^j$ for a suitable $\rho >1.$\\
In order to apply Theorem \ref{FNa} we need an upper bound for
$b_n:=\E\|S'_n - \E S'_n\|.$ Using essentially the same argument as in the proof of Theorem \ref{FNa} we find that this quantity is less than or equal to
$$ \E \|S_n\| + 2\sum_{i=1}^n \E\|X\|I\{\|X\| \ge c_i\}.$$
On account of fact (\ref{fact2}) and condition (\ref{REG}) we have for large enough $j$,
$$b_{n_{j+1}} \le (1 + \delta)\beta_0 c_{n_{j+1}}\le (1 + 2\delta)\beta_0 c_{n_{j}},$$
 provided we have chosen
 $\rho < (1+2\delta)/(1+\delta).$\\ 
From Theorem \ref{FNa} (where we set $\eta =\delta$) and the $c_r$-inequality it now follows for large $j,$
\begin{eqnarray*}
&&\PP\left\{\max_{1 \le n \le n_{j+1}}\|S'_n - \E S'_n\| \ge \{\alpha_0 +\delta +\beta_0(1+\delta)\}(1+ 2\delta)c_{n_j}\right\}\\
&\le& \exp\left(-\frac{(\alpha_0+\delta)^2(1+2\delta)^2 c_{n_j}^2}{(2+\delta)n_{j+1}H(c_{n_{j+1}})}\right) + 8C(\alpha_0+\delta)^{-3} (1+2\delta)^{-3}n_{j+1}\E\|X\|^3I\{\|X\|\le c_{n_{j+1}}\}/c_{n_j}^3\\
&\le& \exp\left(-\frac{(\alpha_0+\delta)^2 c_{n_{j+1}}^2}{2n_{j+1}H(c_{n_{j+1}})}\right) + 8C(\alpha_0+\delta)^{-3}(1 +\delta)^{-3}n_{j+1}\E\|X\|^3I\{\|X\|\le c_{n_{j+1}}\}/c_{n_{j+1}}^3.
\end{eqnarray*}
Recalling relations (\ref{geom}) and (\ref{fact1})  it is easy now to see that (\ref{series}) holds and the proof of the upper bound is complete.
\subsection{The first lower bound}
We now can assume that $\alpha_0 >0,$ but it is possible  that $\alpha_0 = \infty.$\\
It is obviously enough to show that we have for any $0 < \alpha < \alpha_0$ with probability one
\be
\limsup_{n \to \infty}\|S_n\|/c_n \ge \alpha.
\label{lower}
\ee
W. l. o. g. we assume that
\be
\limsup_{n \to \infty} \PP\{\|S_n\| \ge \alpha c_n\} \le 1/2.
\label{stoch}
\ee
Otherwise, we would have $\PP\{\limsup_{n \to \infty}\|S_n\|/c_n \ge \alpha\} \ge \limsup_{n \to \infty} \PP\{\|S_n\| \ge \alpha c_n\} > 1/2$ which implies (\ref{lower}) via the 0-1 law of Hewitt-Savage. \\
We first prove that under the assumptions (\ref{MOM}) and (\ref{stoch}) we have for any sequence $\{n_j\}$ satisfying condition (\ref{sub}),
\be
\sum_{j=1}^{\infty}\PP\{\|S_{n_j}\| \ge \alpha c_{n_j}\} = \infty.
\label{div}
\ee
To that end we choose for any $j$ a functional $f_j \in B_1^*$ so that
$$ \E f^2_j(X)I\{\|X\| \le c_{n_j}\} \ge (1-\epsilon) H(c_{n_j}),$$
where $0 < \epsilon < 1$ will be specified later on.\\
Set for $j,k \ge 1,$
\begin{eqnarray*}
\xi_{j, k} &:=& f_j(X_k)I\{\|X_k\| \le c_{n_j}\},\\
\xi'_{j, k}&:=& \xi_{j, k} - \E \xi_{j, k}.
\end{eqnarray*}
Then it is easy to see that
\be
\PP\{\|S_{n_j}\| \ge \alpha c_{n_j}\} \ge \PP\{\sum_{k=1}^{n_j} \xi_{j,k} \ge
\alpha c_{n_j}\} - n_j \PP\{\|X\| \ge c_{n_j}\}.
\ee
From assumption (\ref{MOM}) it immediately follows that 
$$\sum_{j=1}^{\infty} n_j \PP\{\|X\| \ge c_{n_j}\} < \infty.$$
Moreover, we have $|\E \xi_{j, k}| \le \E \|X\| I\{\|X\| \ge c_{n_j}\}$ which is in view of fact (\ref{fact2}) of order $o(c_{n_j}/n_j).$ Consequently, in order to prove (\ref{div}) it is enough to show that for a suitable $0 < \epsilon  < 1,$
\be
\sum_{j=1}^{\infty} \PP\left\{\sum_{k=1}^{n_j} \xi'_{j,k} \ge
(1+\epsilon)\alpha c_{n_j}\right\} = \infty. \label{f1}
\ee
To estimate these probabilities we employ a non-uniform bound on the rate of convergence in the central limit theorem (see, e.g., Theorem 5.17 on page 168 of Petrov (1995)). We can conclude that
\be
\PP\left\{\sum_{k=1}^{n_j} \xi'_{j,k} \ge
(1+\epsilon)\alpha c_{n_j}\right\} \ge \PP\{\sigma_j \zeta \ge (1+\epsilon)\alpha c_{n_j}/\sqrt{n_j}\} - A\alpha^{-3}(1+\epsilon)^{-3} n_j \E|\xi'_{j,1}|^3 c_{n_j}^{-3},
\ee
where $\zeta$ is a standard normal variable, $\sigma_j^2 = \mathrm{Var}(\xi_{j,1})$ and $A$ is an absolute constant.\\
Noting that $\E|\xi'_{j,1}|^3 \le 8 \E|\xi_{j,1}|^3 \le 8\E\|X\|^3I\{\|X\| \le c_{n_j}\}$, we can infer from fact (\ref{fact1}) that
$$\sum_{j=1}^{\infty} n_j \E|\xi'_{j,1}|^3 c_{n_j}^{-3} < \infty.$$
Therefore, relation (\ref{f1}) and consequently (\ref{div}) follow if we can show that
\be
\sum_{j=1}^{\infty} \PP\{\sigma_j \zeta \ge (1+\epsilon)\alpha c_{n_j}/\sqrt{n_j}\} = \infty.
\label{f2}
\ee
Let $\mathbb{N}_0 = \{j \ge 1: H(c_{n_j}) \le c^2_{n_j}/n_j^2\}.$ Then it is easily checked that for any $\eta >0,$
\be
\sum_{j \in \mathbb{N}_0}\exp\left(-\frac{\eta c_{n_j}^2}{2n_j H(c_{n_j})}\right) < \infty. \label{f3}
\ee
Furthermore, we have for large $j \not \in \mathbb{N}_0,$
\begin{eqnarray*}
\sigma^2_j &=& \E f_j^2(X)I\{\|X\| \le c_{n_j}\} - (\E f_j(X)I\{\|X\| \le c_{n_j}\})^2\\
&=&\E f_j^2(X)I\{\|X\| \le c_{n_j}\} - (\E f_j(X)I\{\|X\| > c_{n_j}\})^2\\
&\ge&(1-\epsilon) H(c_{n_j}) - (\E\|X\|I\{\|X\| \ge c_{n_j}\})^2 \\
&\ge& (1-2\epsilon) H(c_{n_j}).
\end{eqnarray*}
Here we have used that for large $j$, $\E\|X\|I\{\|X\| \ge c_{n_j}\} \le \sqrt{\epsilon} c_{n_j}/n_j$ (see fact (\ref{fact2})).\\
Employing a standard lower bound for the tail probabilites of  normal random variables, we can conclude that for large $j \not \in \mathbb{N}_0,$
$$
\PP\{\sigma_j \zeta \ge (1+\epsilon)\alpha c_{n_j}/\sqrt{n_j}\}
\ge \exp\left(-\frac{(1+\epsilon)^2\alpha^2 c^2_{n_j}}{2n_j (1-3\epsilon)H(c_{n_j})}\right).
$$
Choosing $\epsilon$ so small that $\alpha (1+\epsilon)/\sqrt{1-3\epsilon}) < \alpha_0$ we obtain (\ref{f2}) from relations (\ref{geom}) and (\ref{f3}). This implies relation (\ref{div}).\\
We are now ready to finish the proof by a standard argument.\\
Set $m_k = \sum_{j=1}^k n_j, n \ge 1, $ where $n_j = [(1+\delta^{-2})^j]$ with $0 < \delta < 1/2.$ \\
Note that we then have $n_{j+1}/n_j \ge \delta^{-2}$ and consequently by (\ref{RE}),
 \be c_{n_{j+1}} \ge \delta^{-1} c_{n_j}, j \ge 1.\label{f4}\ee
Likewise it follows that 
$$m_k \le n_k \left(\sum_{i=0}^{k-1} \delta^{2i} \right) \le n_k/(1-\delta^2).$$
Ik $k$ is large enough we  can conclude from (\ref{REG}) that
\be c_{m_k}/c_{n_k} \le (1+\delta)m_k/n_k \le (1-\delta)^{-1}. \label{f5}
\ee
Define for $k \ge 1,$
\begin{eqnarray*}
F_k&:=& \left\{\|S_{m_k} - S_{m_{k-1}}\| \ge \alpha  c_{n_k}\right\},\\
G_k&:=& \left\{\|S_{m_{k-1}}\| \le 2\alpha \delta c_{n_k}\right\}.
\end{eqnarray*}
Note that on account of relations (\ref{f4}) and (\ref{f5}) we have for large $k,$
$$\PP(G_k) \ge \PP \left\{\|S_{m_{k-1}}\| \le 2\alpha  c_{n_{k-1}}\right\}
\ge \PP \left\{\|S_{m_{k-1}}\| \le 2(1-\delta)\alpha  c_{m_{k-1}}  \right\}.$$
Thus (recall (\ref{stoch})) $\PP(G_k) \ge 1/2$ for large $k$.
In view of (\ref{div}) we have $\sum_{k=1}^{\infty} \PP(F_k) = \infty.$ The events $F_k$ and $G_k$ are independent. Thus we can conclude via Lemma
3.4. of Pruitt (1981) that
$$\PP(F_k \cap G_k \mbox{ infinitely often}) =1.$$
We clearly have,
$$ F_k \cap G_k \subset \left\{\|S_{m_k}\| \ge \alpha (1-2\delta)  c_{n_k}\right\}$$
which is due to relation (\ref{f5})
$$ \subset \left\{\|S_{m_k}\| \ge \alpha (1-2\delta)(1-\delta) c_{m_k}\right\}$$
provided that $k$ is large enough.\\
It follows that with probability one,
$$\limsup_{k \to \infty} \|S_{m_k}\|/c_{m_k} \ge  \alpha (1-2\delta)(1-\delta).$$
Since we can choose $\delta$ arbitrarily small, this implies statement (\ref{lower}). \subsection{The second lower bound}
We now assume that $\alpha_0 + \beta_0 < \infty$. If $\alpha_0 = \infty$ the lower bounds follows from part 4.2 and if $\beta_0=\infty$ we can obtain it from Lemma 1.\\

We use essentially the same argument as in Theorem 7 of de Acosta, Kuelbs and Ledoux (1983). There is a small complication: we cannot show for all sequences $\{c_n\}$ satisfying the above conditions that $\E [\sup_n \|S_n\|/c_n] < \infty.$\\ Therefore, we prove this first for  $S'_n = \sum_{i=1}^n X'_i, n \ge 1,$ where the random variables $X'_n$ are defined as in part 4.1, that is
$$X'_n = X_n I\{\|X_n\| \le c_n\}, n \ge 1.$$
From the upper bound part (see 4.1) it follows that we have with probability one,
\be \limsup_{n \to \infty} \|S'_n\| /c_n\le \alpha_0 + \beta_0 < \infty. \label{fa1}\ee
Since $\sup_n \|X'_n\|/c_n \le 1$ we obtain from Corollary 6.12 in Ledoux and Talagrand (1991) that 
\be \E\left[\sup_n \|S'_n\|/c_n\right] < \infty. \label{fa2}\ee
Using the fact that with probability one, $\limsup_{n \to \infty}\|S_n\|/c_n$ is constant,  Fatou's lemma implies that with probability one,
\be \limsup_{n \to \infty}\|S_n\|/c_n=\limsup_{n \to \infty}\|S'_n\|/c_n\ge \limsup_{n \to \infty} \E \|S'_n\|/c_n. \label{fa3}\ee
In view of (\ref{fact2}) we have as $n \to \infty$
$$|\E\|S'_n\| - \E\|S_n\| | \le  \sum_{i=1}^n\E\|X\|I\{\|X\| \ge c_i\}=o(c_n)$$
and we find that with probability one,
$$  \limsup_{n \to \infty} \|S_n\|/c_n \ge \beta_0.$$
This completes the proof of Theorem 5.

\section{Proof of Theorem 3}
We only prove  Theorem 3. The proofs of  the corollaries are exactly as in the real-valued case and they are omitted. Also  Theorems 1 and 2 follow directly from Theorem 5. \\
In view of Theorem 5 and Lemma 1 we only need  to show that
\be \alpha_0 \le \lambda \label{le} \ee
and
\be  \alpha_0 \ge (1-q)^{1/2}\lambda, \; h \in \mathcal{H}_q \label{ge}.\ee
As in the real-valued case we can infer from (\ref{2.9}) that
$$\limsup_{n \to \infty} (LLn)H(a_n/LLn)/h(n)=\lambda^2/2.$$   
The proof of (\ref{ge}) then goes exactly as in the real-valued case and thus it can be omitted as well. In order to prove the corresponding upper bound in the real-valued case, we used another result, namely Theorem 4 in our previous paper, Einmahl and Li (2005). It is possible to extend this result to Banach space valued random variables as well, but there is also a more direct argument for deriving the upper bound (\ref{le}) which we shall give below.\vspace{.3cm}

\nin{\bf Proof of (\ref{le}).}  If $\lambda = \infty$ the upper bound is trivial. Thus we can assume that $\lambda \in [0,\infty).$\\We have to show for any $\alpha > \lambda,$
\be \sum_{n=1}^{\infty} n^{-1}\exp\left(-\frac{\alpha^2 h(n)}{2 H(a_n)}\right) < \infty. \label{u1}\ee
Set $\delta = (\alpha - \lambda)/3.$ Then we clearly have for large enough $n$, 
$$H(a_n/LLn) \le \frac{(\lambda + \delta)^2}{2} \frac{h(n)}{LLn}.$$
Setting $\mathbb{N}_0 =\{n: H(a_n) - H(a_n/LLn) \le \delta \lambda h(n)/LLn\}$ we get for large $n \in \mathbb{N}_0,$
$$ H(a_n) \le \frac{(\lambda + 2\delta)^2}{2} \frac{h(n)}{LLn}$$ and consequently
\be \sum_{n \in \mathbb{N}_0} n^{-1}\exp\left(-\frac{\alpha^2 h(n)}{2 H(a_n)}\right) < \infty.\label{u2}\ee
Further note that we trivially have,
\begin{equation*}
\sum_{n=1}^{\infty} \frac{H(a_n) -H(a_n/LLn)}{a_n^2 LLn}
\le \sum_{n=1}^{\infty} \frac{\E \|X\|^3I\{\|X\| \le a_n\}}{a_n^3} < \infty.
 \end{equation*}
 The latter series is finite because we are assuming $\E \Psi^{-1}(\|X\|) < \infty.$ (See, for instance, Lemma 5(a), Einmahl (1993).)
It follows that
\be \sum_{n \not\in \mathbb{N}_0} \frac{1}{n(LLn)^2} < \infty. \label{u3}\ee
Condition (\ref{2.9}) implies that for large enough $n,$ and $0 < \epsilon <1,$
$$H(a_n) \le (\lambda^2 +1)\frac{a_n^2 LL{a_n}}{\Psi^{-1}(a_n LL{a_n})}
\le C_{\epsilon} \frac{a_n^2 LL{a_n}}{\Psi^{-1}(a_n) (LL{a_n})^{2-\epsilon}} \le C'_{\epsilon} \frac{h(n)}{(LLn)^{1-\epsilon}},$$
where we have used the fact that $\Psi^{-1}$ is  regularly varying at infinity with index 2. ($C_{\epsilon}$ and $C'_{\epsilon}$ are positive constants.)\\
We now can infer from (\ref{u3}) that
\be  \sum_{n \not \in \mathbb{N}_0} n^{-1}\exp\left(-\frac{\alpha^2 h(n)}{2 H(a_n)}\right) \le \sum_{n \not \in \mathbb{N}_0}n^{-1}\exp( -\alpha^2(2C'_{\epsilon})^{-1}(LLn)^{1-\epsilon}) <\infty. \label{u4} \ee
Combining (\ref{u2}) and (\ref{u4}) we see that the series in (\ref{u1}) is finite and our proof of (\ref{le}) is complete. $\Box$\\

\begin{center}
{\bf REFERENCES}
\end{center}
\begin{description}
\item
{\sc Acosta, A. de, Kuelbs, J. and Ledoux, M.} (1983) An inequality for the law of the iterated logarithm. In: Probability in Banach spaces 4.  {\it Lecture Notes in Mathematics} {\bf 990} Springer, Berlin Heidelberg, 1--29.
 \item
{\sc Bousquet, O.} (2002)  A Bennett concentration inequality and its application to suprema of empirical processes.  {\em C. R. Math. Acad. Sci. Paris}  {\bf 334},  no. 6, 495--500. 
 \item
{\sc Einmahl, U.} (1993) Toward a general law of the iterated logarithm
in Banach space. {\em Ann. Probab.} {\bf 21}, 2012-2045.
\item
{\sc Einmahl, U.} (2007) A generalization of Strassen's functional LIL.  {\em J. Theoret. Probab.}, to appear.
\item
{\sc Einmahl, U. and Li, D.} (2005) Some results on two-sided LIL behavior.
\emph {Ann. Probab.}  {\bf 33}, 1601--1624.
\item
{\sc Feller, W.} (1968). An extension of the law of the iterated logarithm to variables without variance. {\em J. Math. Mech.} {\bf 18}, 343-355.
\item
{\sc Klass, M.}  (1976) Toward a universal law of the iterated logarithm  
I.  {\em Z. Wahrsch. Verw. Gebiete} {\bf 36}, 165-178.
\item
{\sc Klass, M.} (1977) Toward a universal law of the iterated logarithm
II.  {\em Z. Wahrsch. Verw. Gebiete} {\bf 39}, 151-165.
\item
{\sc Klein, T. and Rio, E.} (2005) Concentration around the mean for maxima of empirical processes. {\em Ann. Probab.} {\bf 33}, 1060--1077.
\item
{\sc Kuelbs, J.} (1985) The LIL when $X$ is in the domain of attraction of a Gaussian law. {\em Ann. Probab.} {\bf 13}, 825--859.
\item
{\sc Ledoux, M.} (1996) On Talagrand's deviation inequality for product measures. {\em ESAIM Probab. Statist.} {\bf 1}, 63--87.
\item
{\sc Ledoux, M. and Talagrand, M.} (1988) Characterization of the law of the iterated logarithm in Banach space. {\em Ann. Probab.} {\bf 16}, 1242--1264. 
\item
{\sc Ledoux, M. and Talagrand, M.}(1991) {\em Probability in Banach Spaces.} Springer, Berlin (1991).
\item
{\sc Massart, P.} (2000) About the constants in Talagrand's concentration inequalities for empirical processes. {\em Ann. Probab.} {\bf 28}, 863--884.
\item
{\sc Petrov, V. V.} (1995) {\em Limit Theorems of Probability Theory:
Sequences of Independent Random Variables.}~ Clarendon Press,
Oxford. 
\item
{\sc Pruitt, W.} (1981) General one-sided laws of the iterated logarithm.
{\em Ann. Probab.} {\bf 9}, 1--48.
\item 
{\sc Rio, E.:} (2002) Une in\'egalit\'e de Bennet pour les maxima de processus empiriques. {\em Ann. Inst. H. Poincar\'e Probab. Statist.} {\bf 38}, 1053--1057.
\item
{\sc Talagrand, M.:} (1994) Sharper bounds for Gaussian and empirical processes. {\em Ann. Probab.} {\bf 22}, 28--76.
\end{description}

\end{document}